\documentclass[10 pt]{amsart}
\usepackage{amscd,amsmath,amsthm,amssymb}
\usepackage{enumerate,varioref}
\usepackage{epsfig}
\usepackage{graphicx}
\newtheorem{thm}{Theorem}

\newtheorem{prop}[thm]{Proposition}
\theoremstyle{definition}

\theoremstyle{remark}

\newtheorem{rem}[thm]{Remark}

\numberwithin{equation}{subsection}
\newcommand{\R}{\mathbb{R}}

\newcommand{\C}{\mathbb{C}}

\begin{document}

\title{The Second Order Effect of the Quantum Weyl Algebra on a Free Particle }
\author{Clark Alexander}

\email{clark@imsc.res.in} \maketitle \tableofcontents
\section{Introduction}
The goal of the present article is to revisit and add to pair of papers from the mid 1990s which have since received little attention.  The papers [CJSS] and [DaP] both explored the idea of non-relativistic free particles in a slightly noncommutative space.  The former found only the first order effect of deforming a three dimensional space, and only explored the free particle.  The latter further deformed the time coordinate and explored not only a free particle in nonrelativistic quantum mechanics, but also explored the Klein-Gordon equation as well as the Dirac equation in similar situations.  It seems the deformation used for these papers was not popular because it failed to obey a symmetry of $SO_q(3)$ in the case of [CJSS] and $SO_q(3,1)$ in the case of [DaP].  Several different deformations have since replaced these and their use is much more widespread (cf [KS]).  In this article we wish to take the calculations further, to the second order, in which a more peculiar phenomenon occurs.  Furthermore we wish to exhibit the fact that the algebra in question does obey a quantum group symmetry, but one of $Sp_q(6)$ rather than that of $SO_q(3)$.

\section{Preliminaries on the Quantum Weyl Algebra}

The quantum Weyl algebra as it is now so called was first seen explicitly in [WZ].  It is now given as a covariant differential calculus on $\mathcal{O}(\C^N_q)$ named $\mathcal{A}_q(N)$ generated
by $\{X_j,\partial_j | i,j=1,\dots N\}$  with the following relations:

\begin{eqnarray}
X_i X_j &=& q X_j X_i,\hspace{5mm}  i<j,\\
\partial_i \partial_j &=& q^{-1} \partial_j \partial_i,\hspace{5mm}  i<j, \nonumber\\
\partial_i X_j &=& q X_j \partial_i,\hspace{5mm}  i\neq j, \nonumber\\
\partial_i X_i - q^2 X_i \partial_i &=& 1 + (q^2-1)\sum_{j>i}X_j\partial_j. \nonumber 
\end{eqnarray}

\subsection{Quantum Symmetries}

In both [CJSS] and [DaP] similar comments are made on the fact that $\mathcal{A}_q(3)$ does not obey $SO_q(3)$ symmetry.  It requires very little to see that these statements are indeed valid.  Fortunately, the subsets 
$\{ X_i | i=1,\dots, N\}$ and $\{ \partial_i | i=1,\dots,N\}$ obey $GL_q(N)$ and $GL_{q^{-1}}(N)$ symmetries respectively.  If one dives into the literature ever so slightly deeper, one can find the presentation of [KS] in \S 12.3.3 exhibiting 
$\mathcal{A}_q(N)$ as a left quantum space of $Sp_q(2N)$.  

For much of this article we shall only be concerned with $N=3$ and so we will now give a brief synopsis of how $Sp_q(6)$ acts on $\mathcal{A}_q(3)$.

Let us redefine the elements of $\mathcal{A}_q(3)$ as $\{ y_i | i=1,\dots,6\}$ where 
\begin{eqnarray}
y_1 = \alpha q \partial_3, & y_2 = \alpha q^2 \partial_2, & y_3 = \alpha q^3 \partial_1, \\
y_4 = X_1, & y_5 = X_2, & y_6 = X_3, \nonumber
\end{eqnarray}
where $\alpha \in \C^{\times}$.

\begin{prop}
Let $\hat R = (\hat R^{ij}_{kl})$ be the $R$-matrix and $C=(C^i_j)$ the metric for $Sp_q(6)$.  Then the generators $\{X_1,X_2,X_3,\partial_1,\partial_2,\partial_3\}$ satisfy 
the relations for $\mathcal{A}_q(3)$ if and only if
\begin{equation}
\sum_{k,l}^6 \hat R^{ij}_{kl}y_ky_l - q y_iy_j - \alpha q^{-3}C^i_j = 0, \hspace{5mm} i,j=1,\dots,6.
\end{equation}
\end{prop}

The $R$-matrix and metric are given explicitly in [KS] \S 9.3.
The proof of this proposition is essentially a direct calculation and verification that the above  equation reduces to
\begin{equation}
y_{4-j}y_j - q^{-2}y_j y_{4-j} = -q^{-j}\alpha (q^{-2}-1) \sum_{k<j}q^{k-j}y_k y_{4-k}.
\end{equation}

\subsection{A Concrete Realization of $\mathcal{A}_q(3)$}
Since we shall only be working within a three dimensional system, for the remainder of this paper we will adopt the notations that capital letters $(X,Y,Z)$ will represent noncommutative coordinates and lower case letters $(x,y,z)$ will represent commutative coordinates.  Furthermore we will also use $X_1 = X, X_2=Y,X_3=Z$ as well as $x_1=x,x_2=y,x_3=z$ and use the two interchangeably when convenient.

Define the following quantities
\begin{eqnarray}
q &:=& e^{i\theta}, \\
M_j &: =&  x_j \partial_{x_j}, \nonumber\\
M_{k>j} &:=& \sum_{k>j}M_j,\nonumber\\
\beta_j &:=& \left\{ \frac{q^{2(M_j+1)}-1}{(q^2-1)(M_j+1)} \right\}^{1/2}.\nonumber
\end{eqnarray}
Here we restrict to $\theta\in\R$.  We also note that in [CJSS] and [DaP] the operators $M_j$ are written as $N_j$.  However in this paper, as a manner of avoiding confusion, 
we shall use $M_j$ since $N_j$ generally represent number operators in quantum mechanics.

In this scenario we now may write our noncommutative coordinate system explicitly as

\begin{eqnarray}
X_j & = & x_j \beta_j  q^{M_{k>j}},\\
\partial_{X_j} &=& q^{M_{k>j}} \beta_j \partial_{x_j}.
\end{eqnarray}

Expanding this for the sake of clarity we see
\begin{eqnarray}
X = x \beta_1 q^{M_2+M_3},                              &   Y = y \beta_2 q^{M_3},                              & Z = z\beta_3,\\
\partial_X = q^{M_2+M_3} \beta_1 \partial_x, &  \partial_Y = q^{M_3} \beta_2 \partial_y,  & \partial_Z = \beta_3 \partial_z.
\end{eqnarray}

It is a relatively easy computation to verify that these operators satisfy the relations of our algebra $\mathcal{A}_q(3)$.

\section{Slightly Noncommutative Space}
The pair of papers from which the inspiration for this work was derived looked at the algebra $\mathcal{A}_q(3)$ from a physical point of view by considering what happens when $\theta^2 \approx 0$.
In this case the noncommutative parameter becomes $q = 1+ i\theta$.  The commutation relations are made generally simpler, but it creates some interesting physical interpretations of a noncommutative space.
Given that $\theta <<1$ we shall call our new space a \emph{slightly noncommutative space} or \emph{SNCS}.
\begin{rem}
In [DaP] the deformation involved time as well as the three spatial dimensions and therefore the new algebraic system with $q=1+i\theta$ is called a noncommutative space-time (SNCST).
\end{rem} 
In this case we shall look for the second order effect.  For the sake of mathematics, we shall not say $\theta^3 = 0$, but simply ignore the effects of terms involving powers of $\theta^3$.  Including more and more powers should yield the same results when the higher terms are ignored, but the actual calculation is highly more laborious.

Let us begin by looking at the second order expansion of the $\beta_j$.  Recall
\[
\beta_j = \left\{ \frac{    q^{2(M_j+1)} -1 } {(q^2-1)(M_j+1)} \right\}^{1/2} = \left\{ \frac{  q^{2(x_j\partial_j+1)}-1}{(q^2-1)(x_j\partial_j+1)} \right\}^{1/2}.
\]

Tackling the second order expansion of this is best done by expanding the denominator first.  We will have several terms involving powers of theta for which me must compensate in the numerator.  For this reason
it shall appear as if the numerator carries a power of $\theta^3$, however this will be cancelled, and we will be left with a $\theta^2$ as the highest power.

Thus we have in the denominator
\[
(q^2-1)(M_j+1) \approx (1 + 2i\theta + \frac{1}{2}(2i\theta)^2 - 1)(M_j+1) = 2i\theta(M_j+1)(1+i\theta).
\] 

In the numerator
\begin{eqnarray*}
q^{2(M_j+1)}-1 \approx 1 + 2i\theta(M_j+1) + \frac{1}{2}(2i\theta(M_j+1))^2 + \frac{1}{6}(2i\theta(M_j+1))^3 - 1\\ 
= 2i\theta(M_j+1)(1+ i\theta(M_j+1) - \frac{2}{3}\theta^2(M_j+1)^2).
\end{eqnarray*}

Now we cancel the overall factor of $2i\theta (M_j+1)$ and have
\begin{equation}
\beta_j \approx \left\{   \frac{1+ i\theta(M_j+1) - \frac{2}{3}\theta^2(M_j+1)^2}{ 1+i\theta }        \right\}^{1/2}
\end{equation}

From here, we employ the power series representations of 
\[
\frac{1}{1-t} \approx 1 + t + t^2 + O(t^3)
\]

and 

\[
\sqrt{1+t} \approx 1 + \frac{1}{2}t - \frac{1}{8}t^2 + O(t^3)
\]

to arrive at the overall second order approximation of $\beta_j$ as

\begin{equation}
\beta_j \approx  1 - \frac{1}{3}\theta^2 + \frac{1}{2}i\theta x_j\partial_{x_j}  - \frac{3}{8}\theta^2 x_j\partial_{x_j} -\frac{5}{24}\theta^2 x_j^2\partial_{x_j}^2.
\end{equation}

The approximation of $q^{M_{k>j}}$ is nearly trivial.  We have
\[
q^{M_{k>j}} \approx 1+ i\theta M_{k>j} - \frac{1}{2}\theta^2M_{k>j}^2.
\]

To obtain our SNCS coordinates is now simply a matter of multiplication and  bookkeeping. 

\section{The Second Order Effect on a Free Particle}

Up until now we have not explored the idea of making a quantum mechanical system out of our algebra or out SNCS.  The procedure adopted by both [CJSS] and [DaP] is to simply
say we map a Hamiltonian in commutative coordinates into a Hamiltonian in noncommutative coordinates by a simple procedure which looks as follows:

\begin{eqnarray}
x_j &\mapsto & X_j,\\
\partial_{x_j}&\mapsto & \partial_{X_j},\nonumber \\
P_{X_j} & = & -i\hbar \partial_{X_j}, \nonumber\\
H(p_{x_j},x_j) & \mapsto & \mathcal{H}(P_{X_j},X_j).\nonumber
\end{eqnarray}

Given this procedure, we now write the noncommutative Hamiltonian of a free particle as
\begin{equation}
\mathcal{H} = P_X^2 + P_Y^2 + P_Z^2.
\end{equation}

In this new context, the equation we wish to solve is 
\begin{equation}
i\hbar \frac{\partial \Psi}{\partial t} = \frac{1}{2m}(P_X^2+P_Y^2+P_Z^2)\Psi.
\end{equation}

It is not known whether or not this system has a general solution which is unique.  The most natural thing one can attempt 
in this scenario is successive approximations with higher order terms.  The mathematics which follows requires some trickery and 
is at first seemingly unnatural.  In the calculation of the effect of this system on a free particle [CJSS] and [DaP] assume the solution of a slowly varying free particle.
That is 
\begin{equation}
\Psi = \phi(r) e^{i(k\cdot r -\omega t)}
\end{equation}
where $r = (x,y,z)$ the radial vector, $k= (k_1,k_2,k_3)$ the wave vector, and $phi$ a slowly varying function.  That is
\[
| (i/ \phi) \partial_{x_j}\phi | << |k_j|.
\]

The assumption on $\phi$ is a purely physical one.  For mathematical purposes we may approach the calculations as if $\phi=1$.  The assumption of the solution of a free particle is also a physical one.
Perhaps a better physical intuition is assuming that we start in the standard quantum mechanical setting of
\[
i\hbar \partial_t \Psi = \frac{1}{2m}(p_x^2+p_y^2+p_z^2)\Psi
\]

And assuming that we have such an experimental set up in which we may slowly change the physical setting.  In this way we will (in theory) observe the effects on our particle.
In this scenario we shall expand our operator $P_{X_j}$ to the second order and let them hit the assumed solution.  The way in which we collect terms is to leave $p_{x_j}$ alone, and let all the operators hit the wave function
so that we see what extra effects occur.  Given this, one will arrive at the formulae

\begin{eqnarray}
P_X \Psi & \approx & (1-\frac{1}{3}\theta^2)p_x\Psi - \hbar k_1 \theta (\frac{1}{2}k_1x + k_2y + k_3z) \nonumber \\
                 &               & - \hbar k_1\theta^2(-\frac{1}{2}k_1k_2xy -\frac{1}{2}k_1k_3xz -\frac{1}{2}k_2k_3yz )\\
                 &               & -\hbar k_1 \theta^2(-\frac{5}{24}k_1^2x^2   - \frac{1}{2}k_2^2y^2 - \frac{1}{2}k_3^2z^2  ) \nonumber\\
                 &               & - i\hbar k_1 \theta^2(\frac{3}{8}k_1x - \frac{1}{2}k_2y - \frac{1}{2}k_3z) \Psi, \nonumber\\
P_Y \Psi & \approx & (1-\frac{1}{3}\theta^2)p_y\Psi - \hbar k_2 \theta( \frac{1}{2}k_2 y + k_3z) \nonumber\\
                &                & -\hbar k_2 \theta^2 (-\frac{1}{2}k_2k_3 yz -\frac{5}{24}k_2^2y^2 -\frac{1}{2}k_3z^2) \\
                &                & -i\hbar k_2 \theta^2 (\frac{3}{8}k_2y - \frac{1}{2}k_3 z)\Psi, \nonumber\\
P_Z \Psi & \approx & (1-\frac{1}{3}\theta^2)p_z\Psi - \frac{1}{2}\hbar \theta k_3^2 z \nonumber\\
                &                & - \hbar k_3 \theta^2 (-\frac{5}{24} k_3^2 z^2 + i \frac{3}{8} k_3 z)\Psi.                 
\end{eqnarray}

Following the lead of [CJSS] we see that our momentum operators resemble the form
\begin{equation}
P_{X_j} \approx p_{x_j} - A_{x_j}.
\end{equation}

Where $A_{x_j}$ are components of the magnetic potential.  Given this, we assume that our SNCS produces the effect of a magnetic field on our free particle.
Our magnetic field ($B = \nabla \times A$) now appears as 
\begin{eqnarray}
B_x & = & -\hbar \theta k_2k_3  + \hbar\theta^2 k_2( \frac{1}{2}k_2k_3 y + k_3^2 z) + \frac{1}{2}i\hbar\theta^2 k_2k_3,\\
B_y & = & \hbar \theta k_1k_3  - \hbar\theta^2 k_1( \frac{1}{2}k_1k_3 x + \frac{1}{2} k_2k_3 y + k_3^2 z) - \frac{1}{2}i\hbar\theta^2 k_1k_3,  \nonumber\\
B_z & = & - \hbar \theta k_1k_2  + \hbar\theta^2 k_1( \frac{1}{2}k_1k_2 x + \frac{1}{2} k_2k_3 z + k_2^2 y) + \frac{1}{2}i\hbar\theta^2 k_1k_2.\nonumber 
\end{eqnarray}

A couple of interesting phenomena occur in this scenario.  First, is that the first order effect matches exactly with that of [CJSS] and [DaP].  The second is that the new magnetic field is anisotropic.
Since the noncommutative coordinate system `prefers' the $X$ direction in some sense, it seems reasonable that at some point this property of the noncommutative coordinates would be felt by a particle.
Third, is that an imaginary term shows up in our magnetic field.  The physical interpretation of this may be given in several ways:

\begin{enumerate}
\item The imaginary part of $P_{X_j}$ may be pulled out as an imaginary potential term.  This leads to a decay in probability amplitude. \\

\item The imaginary part of the magnetic field is exactly the first order effect with an extra factor of $-i\theta /2$.  In this case we may write
          our magnetic field as 
          \[
          B_i = \hbar \theta \epsilon_{ijk} (1 - \frac{1}{2}i\theta) k_j k_k   + \mathrm{Other\ Second\ Order\ Effects}
          \]  
          Giving rise to the first order magnetic field picking up a phase factor an additional amplitude of $1+ \frac{\theta^2}{4}$.     \\
\item It is also possible that there is simply an imaginary component in the magnetic field.  In Quantum Electrodynamics this can mean that there is a superluminous charge (or faster than light particle), but since we are 
          working in the realm of nonrelativistic quantum mechanics this seems less reasonable.  Therefore, in this case, the physical effect of the imaginary component is more mysterious.                
\end{enumerate}

The main result is that the second order effect of a SNCS on a free particle is the appearance of an position dependent anisotropic field, versus the first order effect which is a simply a constant magnetic field.

\section*{Acknowledgements}
The author wishes to thank Professor R. Jagannathan for his support and openness to questions.  Without his help this work would not have been possible in its present state.  The author also wishes to thank E. Rogers for his help in reading the preliminary work of this article and helping to explain when the physical interpretation was incorrect.


\begin{thebibliography}{10}
\bibitem[CJSS]{CJSS} Chaturvedi, Jagannathan, Sridhar, Srinivasan, \emph{Nonrelativistic Quantum Mechanics in a Noncommutative Space}, J. Phys A: Math. Gen 26 (1993)

\bibitem[DaP]{DaP} Dabrowski, Parashar , \emph{A Free Particle in Noncommutative Space-Time},  Czech J. Phys 46 (1996) 

\bibitem[KS]{KS} Klimyk, Schmudgen, \emph{Quantum Groups and Their Representations}, Texts and Monographs in Physics, Springer-Verlag, 1997

\bibitem[WZ]{WZ} Wess, Zumino, \emph{Covariant Differential Calculus on the $q$-Hyperplane}, Nuclear Physics B, 1990
\end{thebibliography}
\end{document}